\documentclass[11pt]{article}

\usepackage[T1]{fontenc}
\usepackage{lmodern}
\usepackage{microtype}
\usepackage{geometry}
\usepackage{mathtools, amssymb, amsthm}
\usepackage{hyperref}
\hypersetup{colorlinks, allcolors = blue}

\usepackage{enumitem}
\usepackage{xcolor}
\usepackage{graphicx}
\graphicspath{{./comp_resources/}}

\usepackage{caption}
\theoremstyle{plain}

\theoremstyle{definition}

\theoremstyle{remark}

\title{Further verification and empirical evidence for the Erdős-Straus conjecture}

\author{
  Spiridon Mihnea, Bogdan C. Dumitru
  \thanks{Faculty of Mathematics and Computer Science, University of Bucharest.}
  }

\date{August, 2025}

\begin{document}

\maketitle

\begin{abstract}
\noindent We provide empirical evidence for the Erdős-Straus conjecture by improving computational bounds to $10^{18}$ and by evaluating the solution-counting function $f(p)$ for this conjecture.
\end{abstract}

\section{Background}

The Erdős-Straus conjecture states that every fraction of the form $\frac{4}{n}$ can be expanded as the sum of 3 unit fractions $\frac{1}{x}+ \frac{1}{y} + \frac{1}{z}$ with $x,y,z \in \mathbb{N^*}$. The study of this conjecture is concerned with the case where $n$ is a prime number, as unit fraction decompositions for composite numbers $n$ can be obtained from smaller prime numbers: if $n=kp$ for some prime $p$, and Erdős-Straus holds for $p$, then $\frac{4}{p}=\frac{1}{x}+\frac{1}{y}+\frac{1}{z}$ and consequently $\frac{4}{n}=\frac{4}{kp}=\frac{1}{kx}+\frac{1}{ky}+\frac{1}{kz}$.

\section{Extending Salez bounds}

\noindent Many attempts to prove the full conjecture use modular identities involving $p$. For instance, Mordell \cite{mordell} was the first to show that the conjecture is true for all primes, except possibly a small subset given by the modular identity $p \equiv r$ (mod $840$), where $r \in \{1, 121, 169, 289, 361, 529\}$. Currently, the known sets of modular identities are not enough to completely exhaust all possibilities for $p$. However, this approach leads to computational methods that allow the conjecture to be verified up to a large bound. The best-performing method of this type is described by Salez \cite{salez}, whose result we extend.

\subsection{Process}

\noindent Salez defined a modular filter $S_m$ as the set of residue classes mod $m$ for which the conjecture is known to be true and offered an algorithm to produce these filters. Using modular filters, Salez immediately obtains the Mordell result by applying the Chinese remainder theorem on the identities implied by $S_5=\{0,2,3\}$ and $S_7=\{0,3,5,6\}$. By performing this process with the first $7$ prime filters, up to $S_{23}$, Salez obtained the set $R_7$ of residues modulo some $G_7$ that must be checked. Proof for $p\le10^{17}$ follows by verification of integers that escape filtering up to that bound.\\

\noindent We improved this bound to $p\le10^{18}$ by extending this approach with $S_{29}$, obtaining a set $R_8$ with $|R_8|=2101514$ residue classes modulo $G_8=25878772920$ for which we must check the conjecture. Considering the efficiency ratios $\frac{G_7}{|R_7|}$ and $\frac{G_8}{|R_8|}$, this set is roughly twice as efficient.\\

\noindent We divided work in batches $B_k=\{r+kG_8\,|\,r\in R_8\}$ for the sake of multithreading. Verifying the conjecture for all primes $p\le 10^{18}$ is equivalent to checking all batches up to $k = 38641709$, which can be done in parallel. Additionally, the original $10^{17}$ result saves us the need to check the first $k = 3864170$ batches. To verify the integers in any given batch $B_k$, we used Salez' algorithm to precompute a set $\mathcal{S}$ of prime filters with $|\mathcal{S}|=140000$. Then, for each $n\in B_k$, we iterated over each $S_m\in \mathcal{S}$ and checked if $n$ is filtered by $S_m$.

\subsection{Details}

\noindent We note a few things about this process. First, not all integers are filtered by filters in $\mathcal{S}$. We saved these numbers for later processing and found that none of them were prime, therefore they are accounted for by some earlier prime $p$ which was filtered out. Second, some filters are more efficient than others, in that they appear to filter more numbers. We ran our C++ checking program over the first $k=7100$ batches and found that, for instance, $S_{31}$ filtered out a majority of numbers, while most filters were successful $0$ times. After each of these $k$ batches, we sorted the filters according to the total number of integers they filtered. By using the most efficient filters first, we decrease the time it takes to check a batch.\\

\noindent We also remark that computer-aided checking of numbers greater than $10^{17}$ requires us to work around the integer size limits of most programming languages. We generated $R_8$ using a Python rewrite of Salez' algorithm, as the language does not have integer limits, and checked the remaining integers in C++ using the arbitrary-precision integer library GMP \footnote{\url{https://github.com/esc-paper/erdos-straus}}. The inability to use a machine integer for calculations incurred a significant runtime penalty. Our process completed in about 2 weeks with a medium setup.

\section{Solution counting}

Another approach to the Erdős-Straus conjecture is based on a solution-counting function $f(p)=|\{(x, y, z)\,|\,\frac{4}{p} = \frac{1}{x}+\frac{1}{y}+\frac{1}{z}\}|$ for $p\in\mathbb{N^*}$. Elsholtz and Tao \cite{elsholtz-tao} proved that $f(p)$ is upper-bounded polylogarithmically. Furthermore, Bradford \cite{bradford} shows that for any given $p$, all possible $x$ belong to a finite search space $\lceil\frac{p}{4}\rceil\le x\le\lceil\frac{p}{2}\rceil$ and provides an explicit construction of $y$ and $z$ from $x$, given the existence of some divisor $d\mid x^2$ that verifies one of two identities, which we may term the Bradford conditions, depending on the type of the solution, that is, if $p\nmid y$ (Type-1) or if $p\mid y$ (Type-2). This allows us to evaluate $f(p)$, although we remark this is computationally expensive for large $p$.\\

\noindent The Erdős-Straus conjecture itself is equivalent to the statement $f(p)>0\;\,\forall p\in\mathbb{N^*}$. We considered only the ordered set $\mathcal{P}=\{p\;|\;p\,\textrm{prime},\, p \equiv r\, \textrm{(mod 840)},\,r \in \{1, 121, 169, 289, 361, 529\}\}$, as for all other primes the conjecture is known to be true \cite{mordell}. We evaluated $f(\mathcal{P}_i)$ for $i\in\overline{1,N}$, where $N=66737$ and $\mathcal{P}_i$ denotes the $i$-th element of $\mathcal{P}$. This corresponds to the "difficult" primes $p\le3.5\cdot10^7$.\\

\noindent Counting the number of divisors verifying the Bradford conditions across all possible $x$ for some prime $p$ is equivalent to computing $f(p)$. We checked these conditions for a total of $T=29860049601808$ divisors of squares of allowable $x$ for primes $p$ in our considered subset of $\mathcal{P}$ and found that $S=18601583$ of those divisors satisfied at least one of the identities, producing one valid solution to the conjecture for $p$.\\

\begin{figure}[h]
    \centering
    \includegraphics[width=0.8\textwidth]{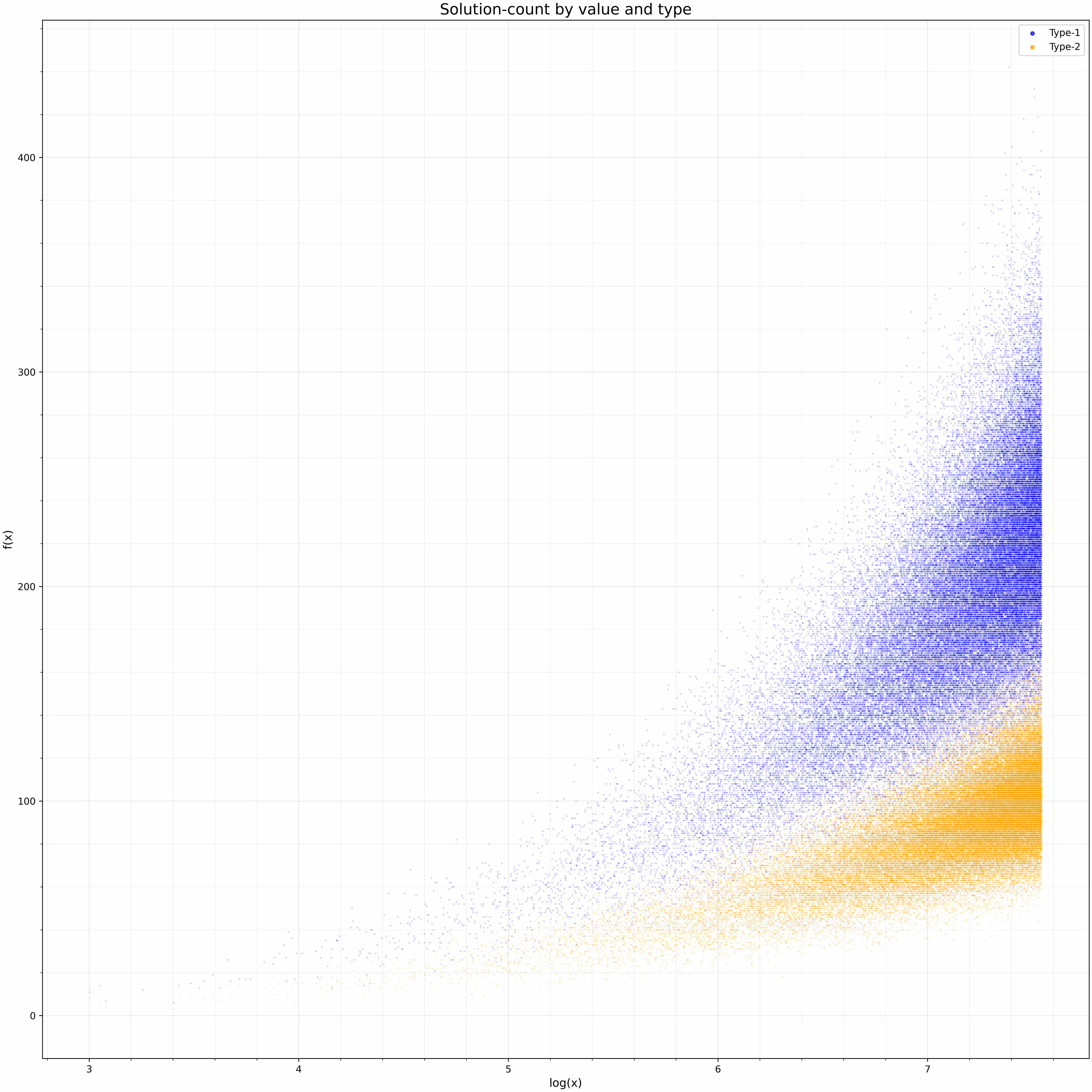}
    \caption{Semi-logarithmic scatter plot of $f(\mathcal{P}_i)$ for $i\in\overline{1,N}$, by solution type}
\end{figure}

\noindent Based on this trial, we find that, empirically, $f(p)$ appears to be increasing consistent with the Elsholtz-Tao upper bound, and furthermore that solutions of Type-1 abound relative to those of Type-2, having found $S_1=12763383$ solutions of Type-1 and only $S_2=5838200$ of Type-2. Figure 1 shows a scatter plot of our $f(p)$ data.

\pagebreak

\section{Bibliography}

\bibliographystyle{plain}
\bibliography{es_comp.bib}

\end{document}